\documentclass[12pt, leqno]{amsart}
\usepackage{amsmath, amsthm, amssymb, amscd}
\usepackage{fancyhdr}
\usepackage{mathdots}
\usepackage{setspace}
\usepackage[small,nohug,heads=vee]{diagrams}
\usepackage{young}
\usepackage{epic}
\diagramstyle[labelstyle=\scriptstyle]
\usepackage[all]{xy}

\theoremstyle{plain}

\newtheorem{thm}{Theorem}[section]

\newtheorem{lem}[thm]{Lemma}
\newtheorem{corl}[thm]{Corollary}
\newtheorem{conj}[thm]{Conjecture}
\theoremstyle{definition}\newtheorem{mydef}[thm]{Definition}
\theoremstyle{remark}\newtheorem{rmk}[thm]{Remark}
\newenvironment{prf}{\textbf{Proof:}}{\rule{1.6ex}{1.6ex}}

\date{\today}

 \newcommand{\Exp}{\operatorname{Exp}}
\newcommand{\Z}{\mathbb{Z}}
\newcommand{\C}{\mathbb{C}}

\newcommand{\g}{\mathfrak{g}}
\newcommand{\h}{\mathfrak{h}}
\newcommand{\hf}{\mathfrak{h}}
\newcommand{\kf}{\mathfrak{k}}
\newcommand{\tf}{\mathfrak{t}}

\newcommand{\slf}{\mathfrak{sl}}

\newcommand{\Lam}{\Lambda}
\newcommand{\lam}{\lambda}

\newcommand{\Del}{\Delta}

\newcommand{\Ad}{\mbox{Ad}}

\newcommand{\la}{\langle}
\newcommand{\ra}{\rangle}

\newcommand{\tl}{\widetilde}

\newcommand{\fb}{\mathfrak{b}}
 \newcommand{\fk}{\mathfrak{k}}
\newcommand{\fg}{\mathfrak{g}}

 \newcommand{\ft}{\mathfrak{t}}

\newcommand{\bc}{\mathbb{C}}

\begin{document}

\title[Reduction of Hitchin's conjecture to the simply-laced case]{Hitchin's conjecture for simply-laced Lie algebras implies that for any simple Lie algebra}

\author{Nathaniel Bushek and Shrawan Kumar}

\address{Department of Mathematics, UNC-Chapel Hill, CB $\#$ 3250, Phillips Hall, Chapel Hill, NC
27599}

\email{bushek$@$unc.edu; shrawan$@$email.unc.edu}

\begin{abstract} Let $\g$ be any simple Lie algebra over $\C$. Recall that there exists an embedding of $\slf_2$ into $\g$, called a principal TDS, passing through a principal nilpotent element of $\g$ and uniquely determined up to conjugation. Moreover, $\wedge (\g^*)^\g$ is freely generated (in the super-graded sense) by primitive elements $\omega_1, \dots, \omega_\ell$, where $\ell$ is the rank of $\g$. N. Hitchin conjectured that for any primitive element $\omega \in \wedge^d (\g^*)^\g$, there exists an irreducible $\slf_2$-submodule $V_\omega \subset \g$ of dimension $d$ such that   $\omega$ is non-zero on the line $\wedge^d (V_\omega)$. We prove that the validity of this conjecture for  simple simply-laced  Lie algebras implies its validity for any simple Lie algebra. 

Let $G$ be a connected, simply-connected, simple, simply-laced algebraic group and let $\sigma$ be a diagram automorphism
of $G$  with fixed subgroup $K$. Then, we show that the restriction map $R(G) \to R(K)$ is surjective, where $R$ denotes the representation ring over $\mathbb{Z}$. As a corollary, we show that the restriction map in the singular cohomology $H^*(G)\to H^*(K)$ is surjective. Our proof of the reduction of Hitchin's conjecture to the simply-laced case relies on this cohomological surjectivity.

\end{abstract}

\maketitle

\section{Introduction}

Let $\fg$ be a finite dimensional simple Lie algebra over the complex numbers $\bc$ with the associated connected, simply-connected complex algebraic group $G$.
Recall that there is a unique (up to conjugation) embedding of $\slf_2$ into $\g$, called a principal TDS, such that the image passes through a principal nilpotent element of $\g$.
Under the adjoint action of a principal TDS, the Lie algebra $\fg$ decomposes as a direct sum of exactly $\ell$ irreducible 
$\slf_2$-submodules
 $V_1, \dots, V_\ell$ of dimensions
  $2m_1+1, \dots, 2m_\ell+1$ respectively, where
$\ell$ is the rank of $\fg$ and $m_1, \dots, m_\ell$ are the exponents
of $\fg$.

Further, the singular cohomology $H^*(G)=H^*(G, \C)$ with complex coefficients is a Hopf algebra. Let $P(\g)
\subset H^*(G)$ be the graded subspace of primitive elements. Then, $P(\g)$ has a basis in degrees 
 $2m_1+1, \dots, 2m_\ell+1$. We identify $H^*(G)$ with $\wedge (\fg^*)^\fg$ and consider $P(\g)$ as a subspace of 
$\wedge (\fg^*)^\fg$.

  Now, N. Hitchin made the following conjecture \cite{Hi}

 \begin{conj}

Let $\g$ be any simple Lie algebra. For any primitive element $\omega \in P_d \subset \wedge^d(\g^*)^\g$, there exists an irreducible subspace $V_\omega\subset \g$ of dimension $d$ with respect to the principal TDS action such that
\[
\omega|_{\wedge^d(V_\omega)} \neq 0.
\]
\end{conj}

  The main motivation for Hitchin behind the above conjecture lies in its connection
  with the study of polyvector fields on the moduli space $M_G(\Sigma)$ of
  semistable principal $G$-bundles on a smooth projective curve $\Sigma$ of any genus $g > 2$.
  Specifically, observe that the cotangent space at a smooth point $E$ of
   $M_G(\Sigma)$ is isomorphic with $H^0(\Sigma, \fg (E)\otimes \Omega)$, where $\fg (E)$
   denotes the associated adjoint bundle and $\Omega$ is the canonical bundle of the
   curve $\Sigma$. Given a biinvariant differential form $\omega$ of degree $k$
  on $G$, i.e., $\omega
  \in \wedge^k(\fg^*)^\fg$, and elements $\Phi_j \in H^0(\Sigma, \fg (E)\otimes \Omega),
  1\leq j\leq k, \omega ( \Phi_1, \dots, \Phi_k)$ defines a skew form with values
  in the line bundle $\Omega^k$.
  Dually, it defines a homomorphism
  $$\Theta_\omega: H^1(\Sigma, \Omega^{1-k}) \to H^0(M_G(\Sigma), \wedge^k\,T),$$
  where $T$ is the tangent bundle of $M_G(\Sigma)$.

  Now, as shown by Hitchin, the validity of the above conjecture would imply that
  the map $\Theta_\omega$ is injective for any invariant form $\omega
  \in \wedge^k(\fg^*)^\fg$ (cf. \cite{Hi}).

 Any simple Lie algebra $\kf$ can be realized as the fixed point subalgebra of a diagram automorphism  of an appropriate simple simply-laced Lie algebra $\g$.  We prove that  the validity of the conjecture for $\g$ implies the validity for $\kf$.
Thus, one needs to verify the conjecture only  for the
simple Lie algebras of types $A, D$  and $E$. Specifically, we have the following result (cf. Theorem \ref{mainh}).
\begin{thm}

If Hitchin's conjecture is valid for any simply-laced simple Lie algebra $\g$, then it is valid for any simple Lie algebra.

More precisely, if Hitchin's conjecture is valid for $\g$ of type $(A_{2\ell -1}; A_{2\ell}; D_4;E_6)$, then it is valid for $\g$ of type $(C_\ell; B_\ell;G_2;F_4)$ respectively.
\end{thm}
The proof relies on constructing a principal TDS in $\kf$ which remains a principal TDS in $\g$. Moreover, we need to use the surjectivity of the space of primitive elements $P(\g)\to P(\kf)$, which allows us to  lift primitive elements  $\omega_d \in \wedge^d (\kf^*)^\kf$ to primitive elements $\tl{\omega}_d \in \wedge^d (\g^*)^\g$. 

Let $K$ be the algebraic subgroup of $G$ with Lie algebra $\kf$, where $\kf$ is the fixed subalgebra under a diagram automorphism of a simple simply-laced Lie algebra $\g$. Our next main result of the paper (cf. Theorem \ref{31}) asserts
the following.

\begin{thm}
The canonical map $\phi: R(G) \rightarrow R(K)$ is surjective, where $R(G)$ denotes the representation ring of $G$ (over $\Z$).

In particular, the canonical   restriction map $\psi: S^\bullet (\g^*)^\g \rightarrow S^\bullet (\kf^*)^\kf$
is surjective. 
\end{thm}
Finally, we use H. Cartan's transgression map and   the surjectivity of $\psi$ to obtain the desired surjectivity of  $\gamma_o: P(\g) \rightarrow P(\kf)$ and thereby the surjectivity of $\gamma:H^*(G)\to H^*(K)$ (cf. Theorem \ref{32}). In our view, the surjectivity of $\phi, \gamma$ and $\gamma_o$ is of independent interest. 

\vskip2ex
 \noindent
{\bf Acknowledgements:}  The second author is grateful to Nigel Hitchin for explaining his conjecture and to Michel Brion for asking the question answered in Theorem \ref{31}. The first author would like to thank Swarnava Mukhopadhyay for many helpful discussions. Both the authors were supported by the NSF grant number DMS-1201310.

\section{Reduction of Hitchin's conjecture to simply-laced Lie algebras}

Let $\g$ be a simple Lie algebra over $\C$ with the associated connected simply-connected complex algebraic group $G$ (with Lie algebra $\g$). 

\begin{mydef}\label{def1}
A Lie algebra embedding $\varphi: \slf_2 \rightarrow \g$ (or its image) is called a principal TDS if $\varphi(X)$ is a principal nilpotent element of $\g$, i.e., $\Ad G\cdot \varphi(X)$ is the open orbit in the nilpotent cone $\mathcal{N}$ of $\g$. 
\end{mydef}

Here, $\slf_2$ is the Lie algebra of traceless $2 \times 2$ matrices over $\C$ with the standard basis
\[
\begin{array}{llll}
X = \begin{pmatrix}
0 & 1 \\
0 & 0 
\end{pmatrix},
&
Y = \begin{pmatrix}
0 & 0 \\
1 & 0 
\end{pmatrix}
\;\;\mbox{and}
&
H= \begin{pmatrix}
1 & 0 \\
0 & -1
\end{pmatrix}.
\end{array}
\]

Let $\varphi': \slf_2 \rightarrow \g$ be another principal TDS. Then, by a result of Kostant \cite{Ko}, Corollary 3.7, $\varphi'$ is conjugate to $\varphi$, i.e., there exists a $g \in G$ such that 
\begin{equation}\label{21}
\varphi' = \Ad g \cdot \varphi.
\end{equation}

Decompose the adjoint representation of $\g$ with respect to a principal TDS $\varphi$ into irreducible components:
\[
\g = V_1 \oplus V_2 \oplus \cdots \oplus V_\ell,
\]
labeling them so that 
\begin{equation}\label{22}
n_1 \leq \cdots \leq n_\ell, \;\;\mbox{where}\;\;n_i = \dim V_i.
\end{equation}
Then, it is known (cf. \cite{Ko}, Corollary 8.7) that 

(a) $\ell = \mbox{rank of $\g$}$.

(b) Each $n_i$ is an odd integer $2m_i + 1$. Moreover, 
$$m_1\leq m_2 \leq \cdots \leq m_\ell$$ are the exponents of $\g$. (The list of exponents for any $\g$ can be found in \cite{Bo}, Planche I - IX.)

(c) Except when $\g$ is of type $D_\ell$ (with $\ell$ even), each $V_i$ is an isotypical component (in particular, uniquely determined) for the principal
TDS $\varphi$, i.e., $m_1 < m_2 < \cdots < m_\ell$. 

When $\g$ is of type $D_\ell$ (with $\ell$ even), the exponents are: 
$$1, 3, 5, \cdots, \ell-3,\ell-1, \ell-1, \ell+1, \cdots , 2\ell-3.$$ Hence, the isotypical component for the highest weight $2\ell - 2$ is a direct sum of two copies of the irreducible module $V_{\slf_2}(2\ell - 2)$ with highest weight $2\ell -2$.

By the identity (\ref{21}), we see that the decomposition of $\g$ with respect to another principal TDS $\varphi'$ looks like
\begin{equation}\label{21new}
\g = \big(\Ad g\cdot V_1\big) \oplus \big(\Ad g\cdot V_2\big) \oplus \cdots \oplus \big(\Ad g\cdot V_\ell\big).
\end{equation}

\begin{mydef}

Recall that the singular cohomology with complex coefficients $H^*(G) = H^*(G, \C)$ is a Hopf algebra, where the product of course comes from the cup product, and the coproduct $\Del: H^*(G) \rightarrow H^*(G) \otimes H^*(G)$ is induced from the multiplication map $\mu : G \times G \rightarrow G$. 
\end{mydef}

Let $P = P(\g) \subset H^*(G)$ be the subspace of primitive elements, i.e.,
\[
P = \{x \in H^*(G) \; |\; \Del(x) = x \otimes 1 + 1 \otimes x \}.
\]
(Observe that $H^*(G)$ does not depend upon the isogeny class of $G$ and hence the notation $P(\g)$ is justified.)

Since $\Del$ is a graded homomorphism, $P \subset H^*(G)$ is a graded linear subspace. It is well-known that, by a result of 
Hopf-Koszul-Samelson,  $P$ is concentrated in odd degrees and, moreover, the canonical map, induced from the product,
\[
\theta: \wedge^\bullet(P) \rightarrow H^*(G)
\]
is a graded algebra isomorphism. In particular, $P$ generates $H^*(G)$ as an algebra over $\C$. 

We can think of $\wedge (\g^*)$ as the algebra of left invariant $\C$-valued forms on a maximal compact subgroup $G_o$ of $G$. By a result of Koszul  (\cite{K}, Th\'eor\`eme 9.2, Chapitre IV), any $\omega \in \wedge(\g^*)^\g$ is a closed form and, moreover, the induced map (identifying $H^*(G_o)$ with the de Rham cohomology $H^*_{dR}(G_o, \C)$)
\[
\eta: \wedge(\g^*)^\g \xrightarrow{\sim} H^*(G_o) \cong H^*(G)
\]
is a graded algebra isomorphism, where the restriction map $H^*(G) \rightarrow H^*(G_o)$ is an isomorphism since $G_o$ is a deformation retract of $G$. 

Via the isomorphism $\eta$, we identify the graded subspace $P \subset H^*(G)$ of primitive elements with a graded subspace (still denoted by) $P \subset \wedge(\g^*)^\g$. 

For any $d \geq 1$, let $P_d$ be the subspace of $P$ of (homogeneous) degree $d$ elements. Then, by \cite{Ko}, Corollary 8.7,
\begin{equation}\label{23}
\dim P_d = \# \{1 \leq i \leq \ell\;|\; n_i = d\},
\end{equation}
where $n_i$'s (given by (\ref{22})) are the dimensions of irreducible components of $\g$ under the principal TDS action. 

In particular, if $\g$ is not of type $D_\ell$ (with $\ell$ even), then 
\begin{equation}\label{24}
\dim P_d \leq 1
\end{equation}
and $P_d$ is of dimension $1$ if and only if $d$ is equal to one of the $n_i's$. If $\g$ is of type $D_\ell$ (with $\ell$ even), 
\begin{equation}\label{25}
\dim P_d \leq 1 \;\;\mbox{if}\;\; d \neq 2\ell -1,\;\;\mbox{and} \;\;\dim P _{2\ell -1}= 2.
\end{equation}

Fix a principal TDS. Hitchin made the following conjecture (cf. \cite{Hi}).

\begin{conj}

Let $\g$ be any simple Lie algebra. For any primitive element $\omega \in P_d \subset \wedge^d(\g^*)^\g$, there exists an irreducible subspace $V_\omega\subset \g$ of dimension $d$ with respect to the principal TDS action such that
\[
\omega|_{\wedge^d(V_\omega)} \neq 0.
\]
\end{conj}

\begin{rmk}
(a) Unless $\g$ is of type $D_\ell$ (with $\ell$ even), given $\omega \in P_d$, there exists a unique irreducible submodule $V$ of dimension $d$ in $\g$ with respect to the principal TDS. Thus,  $V_\omega$ is uniquely determined.

If $\g$ is of type $D_\ell$ (with $\ell$ even), unless $d = 2\ell-1$, given $\omega \in P_d$, there is a unique irreducible submodule $V$ of dimension $d$ in $\g$. Thus, again $V_\omega$ is uniquely determined (for $d \neq 2\ell -1$).

(b) A different choice of principal TDS results in the irreducible submodules being equal to $\Ad g \cdot V$, for some $g \in G$, and some irreducible submodule $V$ for the original principal TDS. But, since we are only considering forms $\omega \in \wedge^d(\g^*)^\g$ (which are, by definition, $\Ad\, G$-invariant), $\omega|_{\wedge^d(\tiny{\Ad}  g \cdot V)} \neq 0$ if and only if $\omega|_{\wedge^d(V)} \neq 0$.
\end{rmk}

Now, we come to  the main result of this section.

\begin{thm}\label{mainh}
If Hitchin's conjecture is valid for any simply-laced simple Lie algebra $\g$, then it is valid for any simple Lie algebra.

More precisely, if Hitchin's conjecture is valid for $\g$ of type $(A_{2\ell -1}; A_{2\ell}; D_4;E_6)$, then it is valid for $\g$ of type $(C_\ell; B_\ell;G_2;F_4)$ respectively.
\end{thm}

\begin{prf}
Let $\kf$ be a non simply-laced simple Lie algebra. Then, there exists a simply-laced simple Lie algebra $\g$ together with a diagram automorphism $\sigma$ (i.e., an automorphism $\sigma$ of $\g$ induced from a diagram automorphism of its Dynkin diagram) such that $\kf$ is the $\sigma$-fixed point $\g^\sigma$ of $\g$. Moreover, given $\kf$, we can choose $\g$ to be of type given in the statement of the theorem. 
(For more details, see Section \ref{folding} on diagram folding.) In particular, we never need to take $\g$ of type $D_\ell$ except $D_4$. 

Choose a Borel subalgebra $\fb$ of $\g$ and a Cartan subalgebra $\tf \subset \fb$ such that they both are stable under $\sigma$. let $\Del = \{\alpha_1 , \dots, \alpha_\ell\} \subset \tf^*$ be the set of simple roots of $\g$, where $\ell$ is the rank of $\g$. Since $\sigma$ keeps $\fb$ 
and $\tf$  stable, $\sigma$ permutes the simple roots. Let $\{\tl{\beta}_1, \dots, \tl{\beta}_{\ell_\kf}\}$ be a set of simple roots taken exactly one simple root from each orbit of $\sigma $ in $\Del$. Then, the fixed subalgebra $\fb_{\kf} := \fb^\sigma$ is a Borel subalgebra of $\kf$, $\tf_\kf : = \tf^\sigma$ is a Cartan subalgebra of $\kf$ and $\{\beta_1, \dots, \beta_{\ell_\kf}\}$ is the set of simple roots of $\kf$, where $\beta_i := \tl{\beta}_i|_{\tf_\kf}$ (cf. \cite{S}). In particular, $\ell_\kf$ is the rank of $\kf$. 

For any $1 \leq n \leq \ell_\kf$, choose a nonzero element $x_n \in \g_{\tl{\beta}_n}$, where $\g_{\tl{\beta}_n}$ is the root space of $\g$ corresponding to the root $\tl{\beta}_n$. Define
\[
y_n = \sum_{i =1}^{\mbox{ord}(\sigma)} \sigma^i(x_n),
\]
where $\mbox{ord}(\sigma)$ is the order of $\sigma$ (which is $2$ except when $\g$ is of type $D_4$ and $\kf$ is of type $G_2$, in which case it is $3$). If $\tl{\beta}_n$ is fixed by $\sigma$, then $\sigma$ acts trivially on $\g_{\tl{\beta}_n}$ (cf. \cite{S}), hence $y_n$ is never zero. Of course, $y_n \in \kf$ and, in fact, $y_n \in \kf_{\beta_n}$. Define the element $y \in \kf$ by 
\[
y = \sum_{n=1}^{\ell_\kf} y_n.
\]

By \cite{Ko}, Theorem 5.3,  $y$ is a principal nilpotent element of $\kf$ and hence there exists a principal TDS in $\kf$:
\[
\begin{array}{lll}
\varphi: \slf_2 \rightarrow \kf & \mbox{such that} & \varphi(X) = y.\\
\end{array}
\]
Moreover, since 
\[
y = \sum_{n=1}^{\ell_\kf}\sum_{i =1}^{\mbox{ord}(\sigma)} \sigma^i(x_n),
\]
again using  \cite{Ko}, Theorem 5.3, we get that $y$ is a principal nilpotent of $\g$ as well. Hence, $\varphi$ is a principal TDS of $\g$ also. Decompose $\g$ under the adjoint action of $\slf_2$ via $\varphi$: 
\[
\g = V_1 \oplus \cdots \oplus V_{\ell_\kf} \oplus V_{\ell_\kf + 1} \oplus \cdots \oplus V_\ell,
\]
where $V_1 \oplus \cdots \oplus V_{\ell_\kf}$ is a decomposition of $\kf$. 

Take a primitive element $\omega_d \in P_d(\kf) \subset \wedge^d(\kf^*)^\kf$, where $P_d(\kf)$ is the space of primitive elements for $\kf$. By (subsequent) Theorem  \ref{32} , the canonical restriction map $\wedge^d(\g^*) \rightarrow \wedge^d(\kf^*)$ induces a surjection
\[
P_d(\g) \rightarrow P_d(\kf), \;\;\mbox{for any $d > 0$}.
\]
Take a preimage $\tl{\omega}_d \in P_d(\g)$ of $\omega_d$. By (\ref{23})-(\ref{24}), there exists a unique irreducible $\slf_2$-submodule $V_{\omega_d}$ of $\kf$ of dimension $d$. Further, by (\ref{23})-(\ref{25}), there exists a unique irreducible $\slf_2$-submodule $V_{\tl{\omega}_d} \subset \g$ of dimension $d$. (For any $\kf$ not of type $G_2$, the uniqueness of $V_{\tl{\omega}_d}$ follows since we have chosen $\g$ not of type $D_\ell$; for $\kf$ of type $G_2$, $P_d(\kf)$ is nonzero if and only if $d=3, 11$ (cf. $\S$\ref{def1}). Again, for these values of $d$, $\dim P_d(D_4) = 1$.)
Hence, $V_{\omega_d} = V_{\tl{\omega}_d}$. Assuming the validity of Hitchin's conjecture for $\g$, we get that ${\tl{\omega}_d}|_{\wedge^d(V_{\tl{\omega}_d})} \neq 0$. Hence, 
\[
{\omega_d}|_{\wedge^d(V_{\omega_d})} = {\tl{\omega}_d}|_{\wedge^d(V_{\tl{\omega}_d})} \neq 0.
\]
This proves the theorem.
\end{prf}

\section{GIT quotient $G//A\MakeLowercase{d} \,G$ and diagram automorphisms}

Let $\g$ be a simple, simply-laced Lie algebra over $\C$ and let $G$ be the connected, simply-connected complex algebraic group with Lie algebra $\g$. Let $\sigma$ be a diagram automorphism of $\g$ and let $\kf = \g^\sigma$ be the fixed subalgebra. Then, $\kf$ is a simple Lie algebra again. Let $K$ be the connected subgroup of $G$ with Lie algebra $\kf$. In fact, $K=G^\sigma$ (cf. \cite{S}). For the connection of the root datum of $K$ with that of $G$, we refer, e.g., to \cite{S}. 

With this notation, we have the following main result of this section.

\begin{thm}\label{31}
The canonical map $\phi: R(G) \rightarrow R(K)$ is surjective, where $R(G)$ denotes the representation ring of $G$ (over $\Z$).

In particular, the canonical map $K// \Ad\; K \rightarrow G // \Ad \;G$, between the GIT quotients, is a closed embedding. 
\end{thm}

Before we come to the proof of the theorem, we need some notational preliminaries on diagram automorphisms and `diagram folding' (i.e., the process of getting $\kf$ from $\g$). As in Section 2, fix a Borel subalgebra $\fb$ and a Cartan subalgebra $\tf \subset \fb$ of $\g$ stable under $\sigma$. Then, $\fb_\kf : = \fb^\sigma$ (resp. $\tf_\kf : = \tf^\sigma$) is a Borel (resp. Cartan) subalgebra of $\kf$. Let $\Del = \{\alpha_1, \dots, \alpha_\ell\} \subset \tf^*$ be the simple roots of $\g$ and let $\{\tl{\beta}_1, \dots, \tl{\beta}_{\ell_\kf}\}$ be a set of simple roots taken exactly one simple root from each orbit of $\sigma$ in $\Del$. Then, $\Del_\kf := \{\beta_1, \dots, \beta_{\ell_\kf}\} \subset \tf_\kf^*$ is the set of simple roots of $\kf$, where $\beta_i : = \tl{\beta}_i|_{\tf_\kf}$. In the following diagrams, we will make a specific choice of indexing convention in each case of diagram folding. 

\subsection{Diagram Folding: Dynkin diagrams of ($\fg, \fk$).} \label{folding}

\hspace{1mm}
\newline
\newline
$\underline{(A_{2n+1}, C_{n+1}):}$

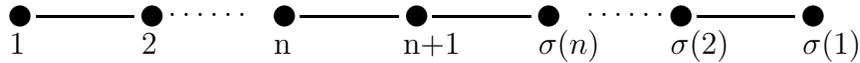
\begin{center}
\begin{picture}(300,24)
\put(0,14){\circle*{8}}
\put(50,14){\circle*{8}}
\put(100,14){\circle*{8}}
\put(150,14){\circle*{8}}
\put(200,14){\circle*{8}}
\put(250,14){\circle*{8}}
\put(300,14){\circle*{8}}

\put(6,14){\line(1,0){37}}
\put(106,14){\line(1,0){37}}
\put(156,14){\line(1,0){37}}
\put(256,14){\line(1,0){37}}

\put(56,14){$\dots\dots$}
\put(214,14){$\dots\dots$}

\put(-4,0){1}
\put(46,0){2}
\put(96,0){n}
\put(145,0){n+1 }
\put(196,0){$\sigma(n)$}
\put(246,0){$\sigma(2)$}
\put(296,0){$\sigma(1)$}


\end{picture}
\end{center}

$\beta_i := {\alpha_i}_{|\ft_\fk}$ for $i \leq n+1$ and $\beta_{n+1}$ is a long root.

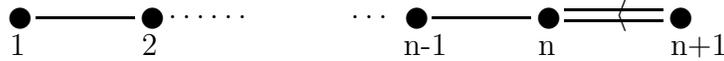
\begin{center}
\begin{picture}(250,24)
\put(0,14){\circle*{8}}
\put(50,14){\circle*{8}}
\put(150,14){\circle*{8}}
\put(200,14){\circle*{8}}
\put(250,14){\circle*{8}}

\put(6,14){\line(1,0){37}}
\put(156,14){\line(1,0){37}}
\put(206,17){\line(1,0){37}}
\put(206,13){\line(1,0){37}}
\put(56,14){$\dots\dots$}
\put(125,14){$\dots$}

\put(-4,0){1}
\put(46,0){2}
\put(145,0){n-1 }
\put(196,0){n}
\put(246,0){n+1}

\put(225, 12){$\la$}
\end{picture}
\end{center}

$\underline{( A_{2n}, B_n):}$

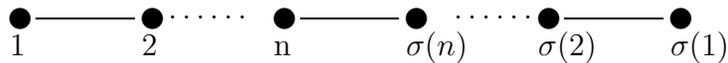
\begin{center}
\begin{picture}(250,25)
\put(0,14){\circle*{8}}
\put(50,14){\circle*{8}}
\put(100,14){\circle*{8}}
\put(150,14){\circle*{8}}
\put(200,14){\circle*{8}}
\put(250,14){\circle*{8}}

\put(6,14){\line(1,0){37}}
\put(106,14){\line(1,0){37}}
\put(206,14){\line(1,0){37}}

\put(56,14){$\dots\dots$}
\put(164,14){$\dots\dots$}

\put(-4,0){1}
\put(46,0){2}
\put(96,0){n}
\put(146,0){$\sigma(n)$}
\put(196,0){$\sigma(2)$}
\put(246,0){$\sigma(1)$}


\end{picture}
\end{center}

$\beta_i = {\alpha_i}_{|\ft_\fk}$ for $1 \leq i \leq n$ and $\beta_{n}$ is a short root.

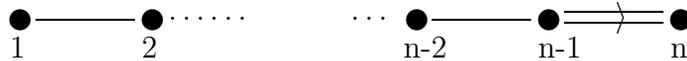
\begin{center}
\begin{picture}(250,24)
\put(0,14){\circle*{8}}
\put(50,14){\circle*{8}}
\put(150,14){\circle*{8}}
\put(200,14){\circle*{8}}
\put(250,14){\circle*{8}}

\put(6,14){\line(1,0){37}}
\put(156,14){\line(1,0){37}}
\put(206,17){\line(1,0){37}}
\put(206,13){\line(1,0){37}}
\put(56,14){$\dots\dots$}
\put(125,14){$\dots$}

\put(-4,0){1}
\put(46,0){2}
\put(145,0){n-2 }
\put(196,0){n-1}
\put(246,0){n}

\put(225, 12){$\ra$}
\end{picture}
\end{center}

$\underline{( D_{n}, B_{n-1}):}$

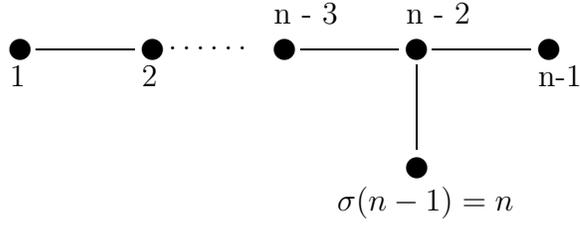
\begin{center}
\begin{picture}(250,82)
\put(0,62){\circle*{8}}
\put(50,62){\circle*{8}}
\put(100,62){\circle*{8}}
\put(150,62){\circle*{8}}
\put(200,62){\circle*{8}}
\put(150,17){\circle*{8}}
\put(6,62){\line(1,0){37}}
\put(106,62){\line(1,0){37}}
\put(156,62){\line(1,0){37}}
\put(150,56){\line(0,-1){32}}

\put(56,62){$\dots\dots$}

\put(-4,48){1}
\put(46,48){2}
\put(96,72){n - 3 }
\put(146,72){n - 2 }
\put(196,48){n-1}
\put(120,1){$\sigma(n-1)=n$}

\end{picture}
\end{center}

$\beta_i := {\alpha_i}_{|\ft_\fk}$ for $1 \leq i \leq n-1$ and $\beta_{n-1}$ is a short root.

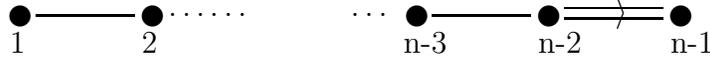
\begin{center}
\begin{picture}(250,24)
\put(0,14){\circle*{8}}
\put(50,14){\circle*{8}}
\put(150,14){\circle*{8}}
\put(200,14){\circle*{8}}
\put(250,14){\circle*{8}}

\put(6,14){\line(1,0){37}}
\put(156,14){\line(1,0){37}}
\put(206,17){\line(1,0){37}}
\put(206,13){\line(1,0){37}}
\put(56,14){$\dots\dots$}
\put(125,14){$\dots$}

\put(-4,0){1}
\put(46,0){2}
\put(145,0){n-3 }
\put(196,0){n-2}
\put(246,0){n-1}

\put(225, 12){$\ra$}
\end{picture}
\end{center}

$\underline{( D_{4}, G_{2}):}$

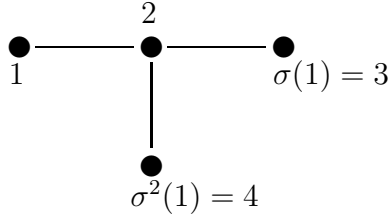
\begin{center}
\begin{picture}(150,80)
\put(0,65){\circle*{8}}
\put(50,65){\circle*{8}}
\put(100,65){\circle*{8}}
\put(50,20){\circle*{8}}

\put(6,65){\line(1,0){37}}
\put(56,65){\line(1,0){37}}
\put(50,59){\line(0,-1){32}}

\put(-4,51){1}
\put(46,74){2}
\put(96,51){$\sigma(1)=3$}
\put(42,5){$\sigma^2(1)=4$}


\end{picture}
\end{center}

$\beta_1 := {\alpha_1}_{|\ft_\fk}$, $\beta_2 := {\alpha_2}_{|\ft_\fk}$, and $\beta_2$ is a long root.

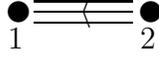
\begin{center}
\begin{picture}(100,20)
\put(0,14){\circle*{8}}
\put(50,14){\circle*{8}}

\put(6,10){\line(1,0){37}}
\put(6,14){\line(1,0){37}}
\put(6,18){\line(1,0){37}}

\put(-4,0){1}
\put(46,0){2}

\put(23, 11){$\la$}
\end{picture}
\end{center}

$\underline{( E_6, F_4):}$

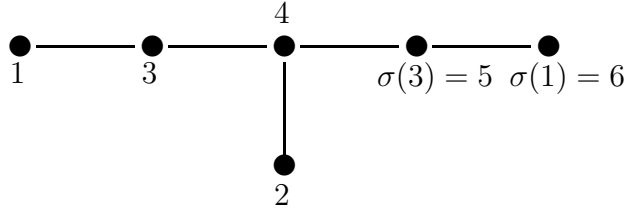
\begin{center}
\begin{picture}(250,80)
\put(0,65){\circle*{8}}
\put(50,65){\circle*{8}}
\put(100,65){\circle*{8}}
\put(150,65){\circle*{8}}
\put(200,65){\circle*{8}}
\put(100,20){\circle*{8}}

\put(6,65){\line(1,0){37}}
\put(56,65){\line(1,0){37}}
\put(106,65){\line(1,0){37}}
\put(156,65){\line(1,0){37}}
\put(100,59){\line(0,-1){35}}

\put(-4,51){1}
\put(46,51){3}
\put(96,74){4}
\put(135,50){$\sigma(3)=5$}
\put(185,50){$\sigma(1) = 6$}
\put(96,5){2}


\end{picture}
\end{center}

$\beta_1 = {\alpha_2}_{|\ft_\fk}$, $\beta_2 = {\alpha_4}_{|\ft_\fk}$, $\beta_3 = {\alpha_3}_{|\ft_\fk}$ and $\beta_4 = {\alpha_1}_{|\ft_\fk}$, with $\beta_2$ a long root.

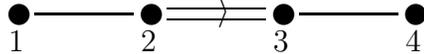
\begin{center}
\begin{picture}(200,20)
\put(0,14){\circle*{8}}
\put(50,14){\circle*{8}}
\put(100,14){\circle*{8}}
\put(150,14){\circle*{8}}

\put(6,14){\line(1,0){37}}
\put(56,12){\line(1,0){37}}
\put(56,16){\line(1,0){37}}
\put(106,14){\line(1,0){37}}

\put(-4,0){1}
\put(46,0){2}
\put(96,0){3 }
\put(146,0){4}

\put(75, 12){$\ra$}
\end{picture}
\end{center}
\vskip2ex

Let $\{\varpi_1, \dots, \varpi_\ell\}$ (resp. $\{\nu_1, \dots, \nu_{\ell_\kf}\}$) be the fundamental weights for the root system of $\g$ (resp. $\kf$). 
We next prove two facts unique to our context. For any simple root $\alpha$, we denote the corresponding coroot by $\alpha^\vee$. We follow the indexing convention as in Subsection \ref{folding}. 

\begin{lem}\label{11}
(a) 
If $G$ is not of type $A_{2n}$ or $E_6$, then $\rho(\varpi_i) = \nu_i$ for $1 \leq i \leq \ell_\fk:=\mbox{rank}(\fk)$. 

(b)
If $G$ is of type $A_{2n}$, then $\rho(\varpi_i) = \rho(\varpi_{2n - i + 1}) = \nu_i$ for $1 \leq i \leq n-1$, and $\rho(\varpi_n)= \rho(\varpi_{n+1}) = 2 \nu_n$. 

(c) If $G$ is of type $E_6$, $\rho(\varpi_1 )= \rho (\varpi_6)  =\nu_4;\,
\rho (\varpi_2) =  \nu_1 ;\,
\rho (\varpi_3) = \rho (\varpi_5 ) = \nu_3; \,
\rho(\varpi_4) = \nu_2 .$
\end{lem}

\begin{prf}
(a)
It suffices to show 
\begin{equation} \label{eqn7} \la \rho(\varpi_i), \beta_j^\vee \ra = \delta_{i,j},\,\,\,\text{ for}\,\,  1\leq i,j\leq \ell_\fk.
\end{equation}
 In this case, we have (\cite{S})
\[
\beta^\vee_j  =
\sum \alpha^\vee_k,
\]
where the summation runs over the orbit of $\alpha_j$ under $\sigma$. For $1\leq j\leq \ell_\fk$, no $\alpha_k$ is in the $\sigma$-orbit of 
$\alpha_j$ for any $1\leq k \leq \ell_\fk$. Thus, the equation \eqref{eqn7} follows. 

(b) 
When $G$ is of type $A_{2n}$, by \cite{S}, 
\[
\beta_j^\vee  =
\begin{cases}
 \alpha_j^\vee + \alpha_{2n - j +1}^\vee ,  &\text{for}\,\, j \leq n-1 ,\\ \\ 
 2\alpha_n^\vee + 2\alpha_{n +1}^\vee,  & \text{for}\,\, j = n.
\end{cases} 
\]
So, for $1 \leq i \leq 2n$, 
\[
\begin{array}{ll}
\la \rho(\varpi_i), \beta_j^\vee \ra
 &=
\begin{cases}
 \la \varpi_i, \alpha_j^\vee \ra + \la \varpi_i, \alpha_{2n - j +1}^\vee\ra, &  \text{for}\,\, j \leq n-1 , \\ \\
2\la \varpi_i, \alpha_n^\vee \ra + 2\la \varpi_i, \alpha_{n +1}^\vee\ra,  & \text{for}\,\,  j = n.
\end{cases}
\\ \\ &
= \begin{cases}
 \delta_{i,j} + \delta_{i , 2n - j +1}, & \text{for}\,\, j \leq n-1 , \\ \\
2 \delta_{i ,n} + 2 \delta_{i , n +1},  & \text{for}\,\, j = n.
\end{cases}
\end{array}
\]
From this (b) follows. 

(c) By \cite{S}, following the indexing convention as in Subsection \ref{folding}, we get that 
\[\beta_1^\vee=\alpha_2^\vee,\,\, \beta_2^\vee=\alpha_4^\vee,\,\, \beta_3^\vee=\alpha_3^\vee + \alpha_5^\vee,\,\, 
\beta_4^\vee=\alpha_1^\vee +\alpha_6^\vee.\]
Thus, \[
\begin{array}{l}
\rho(\varpi_1 )= \rho (\varpi_6)  =\nu_4,\\
\rho (\varpi_2) =  \nu_1 ,\\
\rho (\varpi_3) = \rho (\varpi_5 ) = \nu_3, \\
\rho(\varpi_4) = \nu_2 .\\
\end{array}
\]
\end{prf}

Let $\Lam^+(\g) \subset \tf^*$ (resp. $\Lam^+(\kf) \subset \tf_\kf^*$) be the set of dominant integral weights for the root system of $\g$ (resp. $\kf$) and let $\Lam^+(K) \subset \Lam^+(\kf)$ be the submonoid of dominant characters for the group $K$, i.e., $\Lam^+(K)$ is the set of characters of the maximal torus $T_K$ (with Lie algebra $\tf_\kf$) of $K$ which are dominant with respect to the group $K$. Observe that since $G$ is simply-connected,  $\Lam^+(G) = \Lam^+(\g)$. Moreover, under the restriction  map $\rho: \tf^* \rightarrow \tf_\kf^*$,
\begin{equation}\label{26}
\rho(\Lam^+(\g)) = \Lam^+(K).
\end{equation}

To see this, let $\Lam (K)$ be the character lattice of $K$ (similarly for $\Lam (G) = \Lam (\g)$). Then, by Springer's original construction of $\Lam (K)$ \cite{S}, the restriction $\rho : \Lam (\g) \rightarrow \Lam (K)$ is surjective. Further, from the description of the coroots of $\kf$ as in  \cite{S}, $\rho(\Lam^+(\g)) \subset \Lam^+(\kf)$. Thus, we have $$\rho(\Lam^+(\g)) \subset \Lam^+(\kf) \cap \Lam (K) = \Lam^+(K).$$ 
Conversely, in all cases except for $\g$ of type $A_{2n}$, by Lemma \ref{11},  $\rho(\Lam^+(\g)) = \Lam^+(\kf) \supset \Lam^+(K)$, so equation (\ref{26}) holds in these cases. When $\g$ is of type $A_{2n}$, again by Lemma \ref{11}, 
$$\rho(\Lam^+(\g)) = \bigl(\oplus_{i=1}^{n-1}\,\mathbb{Z}_+ \nu_i\bigr)\oplus 2\mathbb{Z}_+ \nu_n, $$
and
$$\Lambda(K)=\rho(\Lam(\g)) = \bigl(\oplus_{i=1}^{n-1}\,\mathbb{Z} \nu_i\bigr)\oplus 2\mathbb{Z} \nu_n. $$
From this again, we see that \eqref{26} is satisfied. This proves \eqref{26} in all cases.

For any $\lam \in \Lam^+(\g)$, let $V(\lam)$ be the irreducible $G$-module with highest weight $\lam$. Similarly, for $\mu \in \Lam^+(K)$, let $W(\mu)$ be the irreducible $K$-module with highest weight $\mu$. We denote the fundamental representations $V(\varpi_i)$ of $\g$  by $V_i$ and $W(\nu_j)$ of $\kf$  by $W_j$. 

\begin{lem}\label{12}
For any $\lambda \in \Lambda^+(\fg)$,  $W(\rho(\lambda))$ has multiplicity one in $V(\lambda)$
as a $\fk$-module. (Observe that by \eqref{26}, $\rho(\lambda )\in \Lambda^+(K)$.)
\end{lem}
\begin{prf}
Note that the Borel subalgebra $\fb_\fk$ of $\fk$  is contained in the Borel subalgebra $\fb$ of $\g$. So, if $v_\lambda$ is the highest weight vector of $V(\lambda)$ (of weight $\lambda$), then $v_\lambda$ remains a highest weight  vector of weight $\rho (\lambda)$ in $V(\lambda)$ for the action of $\fk$. Hence, $W(\rho(\lambda)) \subset V(\lambda)$.

Multiplicity one is clear from the weight consideration.
\end{prf}

\subsection{Proof of Theorem \ref{31}} Let $\{\mu_1, \dots, \mu_N\} \subset \Lambda^+(K)$ be a set of semigroup generators of 
$\Lambda^+(K)$. Then, the classes $\{[W(\mu_j)]\}_{1\leq j\leq N}$ generate the $\mathbb{Z}$-algebra $R(K)$, where $[W(\mu_j)]
\in R(K)$ denotes the class of the irreducible $K$-module $W(\mu_j)$ (cf. \cite{P}, Theorem 3.12). 

We proceed separately for each of the five cases depending on the type of $(\fg, \fk)$. 

\subsection*{Case I ($A_{2n+1}, C_{n+1}$):}
By Lemmas \ref{11} and \ref{12},  for $1 \leq j \leq n+1$, $W_j \subset  V_j $ (as $\fk$-modules).  
Recall that $V_1 \simeq W_1 \simeq \C^{2n+2}$ (so $W_1 = V_1$) and  $V_j = \wedge^j V_1$ for all $1\leq j \leq 2n+1$. Also, for $2 \leq j \leq n+1$, $W_j$ is given as the kernel of the surjective $\fk$-equivariant contraction map $\wedge^j  W_1 \rightarrow \wedge^{j-2}  W_1$.
Hence, for $2 \leq j \leq n+1$, in $R(\kf)$  (where $R(\kf)$ is the representation ring of $\kf$), by \cite{FH}, Theorem 17.5,
\[ [W_j] + [\wedge^{j-2}  W_1]= [\wedge^j  W_1].
\]
Thus, 
\[\phi([V_1])=[W_1],\,\,\,\text{and} \,\, \phi([V_j])-\phi([V_{j-2}])=[W_j], \,\,\,\text{for}\,\, 2\leq j\leq n+1,\]
where $V_0$ is interpreted as the trivial one dimensional module $\mathbb{C}$. Thus, the class $[W_j]$ of each fundamental representation 
lies in the image of $\phi$, and hence $\phi$ is surjective.

\subsection*{Case II. ($A_{2n}, B_n$):}

By Lemmas \ref{11} and \ref{12}, for $1 \leq j \leq n -1$, $W_j \subset  V_j $  and $W(2 \nu_n) \subset V_n$ (as $\fk$-modules).  
Recall that $V_1 \simeq W_1 \simeq \C^{2n+1}$ (so $W_1 = V_1$), and  $V_j = \wedge^j V_1$ for all $1\leq j\leq 2n$. Also, $W_j = \wedge^j W_1$ for $1 \leq j \leq n-1$ and $W(2\nu_n) = \wedge^n W_1$ (see, e.g., \cite{FH}, Theorem 19.14). Thus, as $\fk$-modules, 
\[
W_j =  V_j,\;\;j \leq n-1;\;\;\;\;\;W(2\nu_n) =  V_n.
\]
Thus,
\[[W_1], \dots, [W_{n-1}], [W(2\nu_n)]\in \,\,\text{Image}\,\, \phi.\]
By Lemma \ref{11} (b) and the identity \eqref{26}, $\Lambda^+(K)$ is generated (as a semigroup) by $\{\nu_1, \dots, \nu_{n-1}, 2\nu_n\}$. Hence,
$\phi$ is surjective in this case.

\subsection*{Case III. $(D_n, B_{n-1}$):}
Recall that $V_1 \simeq \C^{2n}$ and $W_1 \simeq \C^{2n-1}$. By Lemmas \ref{11} and \ref{12}, for $1 \leq j \leq n -1$, $W_j \subset  V_j $  (as $\fk$-modules).  Since $W_1 \subset  V_1$ (as $\fk$-modules), we get (as $\fk$-modules):
\[
 V_1 = W_1 \oplus \C.
\] 
Thus, for $1 \leq k \leq n-2$, as $\fk$-modules, 
\[ V_k = \wedge^k  V_1  = \wedge^k  (W_1 \oplus \C) \simeq  (\wedge^{k}W_1) \oplus (\wedge^{k-1}W_1) = W_k \oplus W_{k-1},
\]
where the first equality is by \cite{FH}, Theorem 19.2; $W_0$ is interpreted as the one dimensional trivial module and the last equality is from the proof of Case II. 

Since $W_{n-1} \subset V_{n-1}$ as $\fk$-modules, and both being spin representations have the same dimension 
$2^{n-1}$ (see, e.g., \cite{GW}, Section 6.2.2), we get $V_{n-1} = W_{n-1}$. Therefore,
\[\phi([V_k])=[W_k]+[W_{k-1}]\,\,\,\text{for}\,\,1\leq k\leq n-2, \,\,\,\text{and}\,\,\phi([V_{n-1}])=[W_{n-1}].\]
In particular, each of $[W_1], \dots, [W_{n-1}]$ lies in the image of $\phi$, proving the surjectivity of $\phi$ in this case.

\subsection*{Case IV. ($D_4, G_2$):}
The two fundamental representations $W_1$ and $W_2$ have respective dimensions $7$ and $14$ (\cite{FH}, Section 22.3).  On the other hand, $V_1 $ is eight dimensional and $V_2 = \wedge^2 V_1$. 
Since $\rho(\varpi_1) = \nu_1$ (by Lemma \ref{11}),  by Lemma \ref{12} we get  $W_1 \subset  V_1$ (as $\fk$-modules). So, we have the decomposition (as $\fk$-modules):
\[
V_1 = W_1 \oplus \C.
\]
Thus, as $\fk$-modules, 
\[
 V_2 = \wedge^2  V_1 = \wedge^2(W_1 \oplus \C) \simeq \big(\wedge^2 W_1\big) \oplus W_1. 
\]
But, $\wedge^2 W_1 \simeq W_2 \oplus W_1$ (\cite{FH}, Section 22.3). Hence, as $\fk$-modules,
\[
 V_2 = W_2 \oplus W_1^{\oplus 2}.
\] 
This gives
\[\phi([V_1])=[W_1]+1\,\,\,\text{and}\,\,\phi([V_2])=[W_2]+2 [W_1],
\]
which proves the surjectivity of $\phi$ in this case.

\subsection*{Case V. ($E_6, F_4$):} 
By Lemma \ref{11}(c), we see that $\rho$ is surjective with kernel given by $\{a \varpi_1 + b \varpi_3 - b\varpi_5 - a \varpi_6\;|\; a, b \in \Z\}$. 
Considering the images of $\varpi_i$ under $\rho$, we have as $\fk$-modules (by Lemmas \ref{11}(c) and \ref{12}), 
\[
\begin{array}{l}
W_1 \subset  V_2 ,\\
W_2 \subset  V_4 ,\\
W_3 \subset  V_3,  V_5, \\
W_4 \subset  V_1,  V_6.
\end{array}
\]
Using  \cite{Sl}, Tables 44 and 47 or \cite{Li}, we obtain
\[
\begin{array}{ll}
\dim(W_1) = 52, & \dim(V_2) = 78, \\ 
\dim(W_2) = 1274, & \dim(V_4) = 2925, \\ 
\dim(W_3) = 273, & \dim(V_3) = \dim(V_5) = 351, \\ 
\dim(W_4) = 26, & \dim(V_1) = \dim(V_6) = 27.
\end{array}
\]
Along with the fundamental $\fk$-modules, there are only three other irreducible $\fk$-modules of dimensions at most 1651 (\cite{Sl}, Table 44, or \cite{Li}). These are
$\dim(W(2\nu_4)) = 324$, $\dim(W(\nu_1 + \nu_4)) = 1053$, and $\dim(W(2\nu_1)) = 1053$.

Let $U^k$ denote an arbitrary $\fk$-module of dimension $k$. Considering the dimensions, we get (as $\fk$-modules):
\[
\begin{array}{l}
V_1 =  V_6 = W_4 \oplus \C,\\ 
  V_2  = W_1 \oplus U^{26},\\ 
  V_3 = V_5 = W_3 \oplus U^{78},\\ 
 V_4 = W_2 \oplus U^{1651}.
\end{array}
\]
Now, $U^{26}$ must be either $W_4$ or the trivial module $\C^{26}$, and $U^{78}$ must be some combination of $W_4$, $W_1$ and $\C$. 
Since $\phi ([ V_1]) - 1 = [W_4]$, this implies that $[W_4]$, $[W_1]$ and $[W_3]$ are in the image of $\phi$. 
(We remark that \cite{Sl} gives $F_4 \subset E_6$ branching, but we continue without these results for clarity and completeness.)

Using appropriate tensor product decompositions in \cite{Li}, we get 
\begin{align}
[W(2\nu_4)] &=  [W_4]^2 - [W_3] - [W_1] - [W_4] -1, \\  
 [W(\nu_1 + \nu_4) ] &= [W_1] [W_4] - [W_3] - [W_4], \\
[W(2\nu_1)] &= [W_1]^2 - [W_2] - [W(2\nu_4)] - [W_1] - 1. 
\end{align}
Since  $W_2$ appears in $V_4$ as  a $\fk$-submodule exactly once by Lemma \ref{12}, from the above identities, we get that  $[W_2]$ lies in the image of  $\phi$  if  $W(2\nu_1)$ is not a component of  $V_4$. In fact, we prove below that $2 \nu_1$ is not a $\fk$-weight of $V_4$ at all.

In order that $2 \nu_1$ be a $\fk$-weight of $V_4$, we should have
$2 \nu_1 = \mu|_{\ft_\fk}$, where $\mu $ is a weight of $V_4$. This is only possible if there exists a weight of $V_4$ of the form $\mu = a\varpi_1 + 2 \varpi_2 + b \varpi_3 - b \varpi_5 - a \varpi_6$, for some $a,b\in \mathbb{Z}$. We claim this is impossible. Indeed, all weights of $V_4$ are of the form $\varpi_4 - \sum_{i=1}^6 d_i \alpha_i$, where $d_i\in \Z^+$. If such $\mu$ existed, then by \cite{Bo}, Planche V, 
\[
\begin{array}{l}
\sum_{i=1}^6 d_i \alpha_i = \varpi_4 - \mu \\ 
= \varpi_4 + a ( \varpi_6 - \varpi_1) - 2 \varpi_2 + b  (\varpi_5 - \varpi_3) \\ 
=(2 \alpha_1 + 3 \alpha_2 + 4 \alpha_3 + 6 \alpha_4+ 4 \alpha_5 + 2 \alpha_6) 
+(a/3)(-2 \alpha_1 - \alpha_3 + \alpha_5 + 2 \alpha_6)\\ 
-2(\alpha_1 + 2 \alpha_2 + 2 \alpha_3 + 3 \alpha_4 + 2 \alpha_5 + \alpha_6)
+(b/3)(- \alpha_1 - 2 \alpha_3 + 2 \alpha_5 + \alpha_6),
\end{array}
\]
from which we immediately see a contradiction since the $\alpha_2$ coefficient is $-1$.

This completes the proof  in this last case and hence the proof of the first part of Theorem \ref{31} is completed.

To prove that $\eta: K// \Ad\; K \rightarrow G // \Ad \;G$ is a closed embedding, it suffices to show that the induced map between the affine coordinate rings $\eta^*:\mathbb{C}[G // \Ad \;G]\to \mathbb{C}[K // \Ad \;K]$ is surjective. But, by \cite{P}, Theorem 3.5, there is a functorial isomorphism
$$\mathbb{C}\otimes_\mathbb{Z}\,R(G) \to \mathbb{C}[G // \Ad \;G],$$
and  similarly we have an isomorphism
$$\mathbb{C}\otimes_\mathbb{Z}\,R(K) \to \mathbb{C}[K // \Ad \;K].$$
From this the surjectivity of $\eta^*$ follows from the surjectivity of $R(G)\to R(K)$. This proves the theorem.
\qed

We give the following Lie algebra analogue as a corollary.

\begin{corl}\label{35}
The canonical restriction map
$$ S(\fg^*)^\fg \to S(\fk^*)^\fk $$
is surjective.
\end{corl}
\begin{prf} By \cite{St}, $\S$6.4, for any connected semisimple algebraic group $H$ over $\C$, the restriction map
\begin{equation}\label{260}
r:\C[H// \Ad\; H]\simeq \C[H]^H\to \C[T_H]^{W_H}
\end{equation}
is an isomorphism of $\mathbb{C}$-algebras, where $T_H\subset H$ is a maximal torus and $W_H$ is the Weyl group of $H$. 

Similarly, the restriction map 
\begin{equation}\label{261}
r_o:\C[\h]^H\to \C[\tf_{\h}]^{W_H}
\end{equation}
is a graded algebra isomorphism, where $\h$ (resp. $\tf_\h$) is the Lie algebra of $H$ (resp. $T_H$). Thus, to prove the corollary, it suffices to show that the canonical restriction map
\[
\beta_o^*:\C[\tf]^W \rightarrow \C[\tf_\kf]^{W_K}
\] is surjective,where $W$ (resp. $W_K$) is the Weyl group of $G$ (resp. $K$). Since $\beta_o^*$ is a graded algebra homomorphism induced from the 
$\C^*$-equivariant map $\beta_o: \tf_\kf/{W_K} \to \tf/W$ (where the $\C^*$-action is the standard homothety action), it suffices to show that the tangent map between  the Zariski tangent spaces at $0$:
$$(d\beta_o)_0:T_0(\tf_\kf/W_K) \rightarrow T_0(\tf/W)$$
is injective. Let $T^{anal}$ denote the analytic tangent space. Then, the canonical map 
$$T_x^{anal}(X) \to T_x(X)$$
is an isomorphism for any algebraic variety $X$ and any point $x\in X$. 

Consider the commutative diagram:
\[
\begin{diagram}
\tf_\kf/W_K &\rTo^{\hspace{3mm}\beta_o\hspace{3mm}}  & \tf/W\\
\dTo_{\Exp}& & \dTo_{\Exp}\\
T_K/W_K &\rTo_{\hspace{3mm}\beta\hspace{3mm}} & T/W,
\end{diagram}
\]
where $T_K\subset K$ is the maximal torus with Lie algebra $\tf_\kf$ and $\beta:T_K/W_K\to T/W$ is the canonical map. Since $T_K, T$ 
are tori, $\Exp$ is a local isomorphism in the analytic category. In particular, there exist open subsets (in the analytic topology) $0\in U_\kf\subset 
\tf_\kf/W_K, 0\in U\subset \tf/W, 1\in V_K\subset T_K/W_K$ and 
$1\in V\subset T/W$ such that $\beta_o(U_K)\subset U$ and $\Exp_{|U_\kf}:U_\kf\to V_K$ is an analytic isomorphism and so is  
$\Exp_{|U}:U\to V$. Since, by Theorem \ref{31} and the isomorphism \eqref{260}, $\beta$ is a closed embedding, 
$$(d\beta)_1:T_1^{anal}(T_K/W_K)\simeq T_1(T_K/W_K) \to T_1^{anal}(T/W)\simeq T_1(T/W)$$
is injective and hence so is $ T_0(\tf_\kf/W_K)\to T_0(\tf/W)$. This proves the corollary.
\end{prf}

As a consequence  of Corollary \ref{35}, we get the following.

\begin{thm}\label{32}
With the notation and assumptions as in Theorem \ref{31}, the canonical restriction map $\gamma : H^*(G) \rightarrow H^*(K)$ is surjective. Moreover, this induces a surjective (graded) map 
\[
\gamma_o : P(\g) \rightarrow P(\kf),
\]
where $P(\g) \subset H^*(G)$ is the subspace of primitive elements. 
\end{thm}
\begin{prf}
From the definition of coproduct, it is easy to see that the following diagram is commutative:
\[
\begin{diagram}
H^*(G) &\rTo^{\hspace{3mm}\Del_G\hspace{3mm}} & H^*(G) \otimes H^*(G)\\
\dTo_{\gamma}& & \dTo^{\gamma \otimes \gamma}\\
H^*(K) &\rTo_{\hspace{3mm}\Del_K\hspace{3mm}} & H^*(K) \otimes H^*(K).
\end{diagram}
\]
Thus, $\gamma$ takes $P(\g)$ to $P(\kf)$. 

Let $\hf$ be a reductive Lie algebra. For any $v\in \hf$, define the derivation $i(v):S(\hf^*)\to S(\hf^*)$ given by $i(v)(f)=f(v)$, for $f\in \hf^*$. Further, define an algebra homomorphism $\lambda:S(\hf^*)\to \wedge^{\text{even}}(\hf^*)$ by $\lambda(f)=df$, for $f\in \hf^*=S^1(\hf^*)$, where
$d:\wedge^1(\hf^*)=\hf^*\to \wedge^2(\hf^*)$ is the standard differential in the Lie algebra cochain complex $\wedge^\bullet(\hf^*)$. Now, define the {\it transgression map}
$$\tau = \tau_\hf : S^+(\hf^*)^\hf\to \wedge^+(\hf^*)^\hf, \,\,\,\tau(p)=\sum_j\,e_j^*\wedge \lambda(i(e_j)p), $$
for $p\in S^+(\hf^*)^\hf$, where $\{e_j\}$ is a basis of $\hf$ and $\{e_j^*\}$ is the dual basis of $\hf^*$.

By a result of Cartan (cf. \cite{Ca}, Th\'eor\`eme 2; also see [L]), $\tau$ factors through  $$S^+(\hf^*)^\hf/(S^+(\hf^*)^\hf)\cdot (S^+(\hf^*)^\hf)$$ to give an injective map 
$$\bar{\tau}: S^+(\hf^*)^\hf/(S^+(\hf^*)^\hf)\cdot (S^+(\hf^*)^\hf) \to \wedge^+(\hf^*)^\hf$$
with image precisely equal to the space of primitive elements $P(\hf)$. From the definition of $\tau$, it is easy to see that the following diagram is commutative:
\[
\begin{diagram}
S^+ (\g^*)^\g&\rTo^{\hspace{3mm}\tau_\g\hspace{3mm}} & \wedge^+(\g^*)^\g\\
\dTo& & \dTo\\
 S^+ (\kf^*)^\kf &\rTo^{\hspace{3mm}\tau_\kf\hspace{3mm}} &  \wedge^+(\kf^*)^\kf,
\end{diagram}
\]
where the vertical maps are the canonical restriction maps. By using Corollary \ref{35}, this proves that $P(\g)$ surjects onto $P(\kf)$. Since
$P(\kf)$ generates  $\wedge^*(\kf^*)^\kf\simeq H^*(K)$  as an algebra, we get that $\gamma$ is surjective. This proves the theorem.
\end{prf}

\begin{rmk}

As a consequence of the above theorem, we see that the Leray-Serre homology (or cohomology) spectral sequence with coefficients in $\C$ for the fibration 
\[
K \rightarrow G \rightarrow G/K
\]
degenerates at the $E^2$-term.
\end{rmk}

\bibliographystyle{plain}
\def\noopsort#1{}

\end{document}